\title{Maps admitting trialities but not dualities}

\author{Gareth A. Jones\\
School of Mathematics\\
University of Southampton\\
Southampton SO17  1BJ, U.K.\\
{\tt G.A.Jones@maths.soton.ac.uk}
\and
Andrew Poulton\\
School of Mathematics\\
University of Southampton\\
Southampton SO17  1BJ, U.K.\\
{\tt awdp106@soton.ac.uk}
}

\documentclass[12pt]{article}
\usepackage[latin1]{inputenc}
\usepackage{a4wide}
\usepackage{latexsym}
\usepackage{amsmath}
\usepackage{amssymb}
\usepackage{amsfonts}
\newtheorem{thm}{Theorem}[section]
\newtheorem{lemma}[thm]{Lemma}
\newtheorem{cor}[thm]{Corollary}
\newtheorem{prop}[thm]{Proposition}
\date{}
\begin{document} 

\maketitle

\begin{abstract}
\noindent We use group theory to construct infinite families of maps on surfaces which are invariant under Wilson's map operations of order $3$ but not under the operations of order $2$, such as duality and Petrie duality.
\end{abstract}



\noindent{\bf Running head:} Maps admitting trialities.\\

\section{Introduction}

In~\cite{Wil}, Wilson introduced a group $\Sigma$ of six operations on regular maps, isomorphic to the symmetric group $S_3$. In addition to the identity operation, $\Sigma$ contains the duality operation $D$, transposing vertices and faces, and the Petrie duality $P$, transposing faces and Petrie polygons (closed zig-zag paths); these two operations have order $2$, and they generate the group $\Sigma$. The remaining elements of $\Sigma$ are the `opposite' operation $DPD=PDP$, which is a duality transposing vertices and Petrie polygons, and two mutually inverse operations $DP$ and $PD$ of order $3$, called {\sl trialities}, each permuting vertices, faces and Petrie polygons in a cycle of length $3$.

The images of a map under these six operations are called its {\sl direct derivates}, and Wilson divided regular maps into classes I -- IV as they have $6$, $3$, $2$ or $1$ of these, up to isomorphism (the only possibilities). Maps $\cal M$ in class~III are invariant under the triality operations $DP$ and $PD$, but not under the duality operations $D$, $P$ and $DPD$; thus
\[{\cal M}\cong DP({\cal M})\cong PD({\cal M})\quad {\rm and}\quad
D({\cal M})\cong P({\cal M})\cong DPD({\cal M}),\]
but ${\cal M}\not\cong D({\cal M})$. Wilson pointed out in~\cite{Wil} that such maps seem to be very rare, and indeed at one stage he had suspected that they did not exist. The single example given in~\cite{Wil} is a dual pair $\cal W$ and $D({\cal W})$ of non-orientable maps of type $\{9,9\}_9$ and characteristic $-70$, with the simple group $L_2(8)$ of order $504$ as their common automorphism group. These appear to be the only examples currently in the literature, although after Wilson raised the question of their existence at the SIGMAC 06 conference in Aveiro, Conder~\cite{Con06} found several further examples by using a computer program.

In~\cite{JT}, Jones and Thornton showed that regular maps correspond to normal subgroups of a certain group $\Gamma$, the free product $V_4*C_2$ of a Klein four-group $V_4$ and a cyclic group $C_2$ of order $2$, and that the operations in $\Sigma$ are induced by the outer automorphism group ${\rm Out}\,\Gamma\cong {\rm Aut}\,V_4\cong S_3$ acting on such subgroups. Our aim here is to use this algebraic approach to reinterpret Wilson's example, and to give several constructions which yield infinite series of regular maps in class~III.

We thank Marston Conder and Steve Wilson for some valuable comments, and the Nuffield Foundation for supporting the second author with a Nuffield Science Bursary.

\section{Maps}

Here we outline the algebraic theory of maps. For simplicity of exposition we will restrict attention to maps without boundary and without free edges, since it is easy to show that there are no regular maps in class~III possessing these features. For a more comprehensive account, see~\cite{BS}.

If $\cal M$ is a map then let $\Omega$ be the set of {\sl blades\/} associated with $\cal M$: these are local flags $\alpha=(v,e,f)$ where $v$ is a vertex, $e$ is an edge, and $f$ is a face, all mutually incident. We define three permutations $r_0, r_1$ and $r_2$ of $\Omega$ as follows: $r_i$ sends each flag $\alpha$ to an adjacent flag $\alpha r_i$ by preserving its $j$-dimensional components for each $j\neq i$, while changing (in the only possible way) its $i$-dimensional component. These permutations are illustrated in Figure~1, where the small triangles represent blades.

\begin{picture}(100,160)(0,10)

\put(120,60){\circle*{5}}
\put(280,60){\circle*{5}}
\put(110,50){$v$}
\put(195,120){$f$}
\put(195,65){$e$}

\put (120,59.7){\line(1,0){160}}
\put (120,59.8){\line(1,0){160}}
\put (120,59.9){\line(1,0){160}}
\put (120,60){\line(1,0){160}}
\put (120,60.1){\line(1,0){160}}
\put (120,60.2){\line(1,0){160}}
\put (120,60.3){\line(1,0){160}}

\put (120,59.7){\line(-1,2){40}}
\put (120,59.8){\line(-1,2){40}}
\put (120,59.9){\line(-1,2){40}}
\put (120,60){\line(-1,2){40}}
\put (120,60.1){\line(-1,2){40}}
\put (120,60.2){\line(-1,2){40}}
\put (120,60.3){\line(-1,2){40}}

\put (120,60){\line(2,1){30}}
\put (120,60){\line(2,-1){30}}
\put (280,60){\line(-2,1){30}}
\put (280,60){\line(-2,-1){30}}
\put(150,45){\line(0,1){30}}
\put(250,45){\line(0,1){30}}

\put (120,60){\line(0,1){35}}
\put (105,88){\line(2,1){15}}

\put(145,80){$\alpha$}
\put(242,80){$\alpha r_0$}
\put(142,35){$\alpha r_2$}
\put(217,35){$\alpha r_0r_2=\alpha r_2r_0$}
\put(115,100){$\alpha r_1$}

\end{picture}

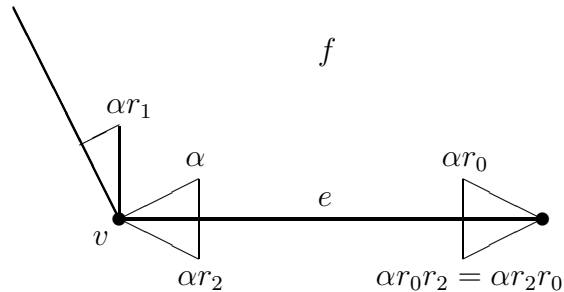
\centerline{Figure 1: the permutations $r_0$, $r_1$ and $r_2$.}

\bigskip

By connectivity these permutations generate a transitive subgroup $G$ of the symmetric group ${\rm Sym}\,\Omega$ on $\Omega$, known as the {\sl monodromy group\/} ${\rm Mon}\,{\cal M}$ of $\cal M$. They satisfy
\[r_0^2=r_1^2=r_2^2=(r_0r_2)^2=1,\]
so that if we define an abstract group $\Gamma$ by means of the presentation
\[\Gamma=\langle R_0, R_1, R_2\mid R_0^2=R_1^2=R_2^2=(R_0R_2)^2=1\rangle\]
then there is an epimorphism $\theta:\Gamma\to G$ given by $R_i\mapsto r_i$. (The generators $R_i$ of $\Gamma$ were denoted by $l, r$ and $t$ in~\cite{JT}, but here we follow the notation of Coxeter and Moser~\cite{CM}.) The subgroups $M=\theta^{-1}(G_{\alpha})$ of $\Gamma$ fixing a blade $\alpha\in\Omega$, known as the {\sl map subgroups}, form a conjugacy class of subgroups of index $|\Omega|$ in $\Gamma$. One can reverse this process: any transitive permutation representation of $\Gamma$ gives rise to a map in which the vertices, edges and faces correspond to the orbits of the dihedral subgroups $\langle R_1, R_2\rangle$, $\langle R_0, R_2\rangle$ and $\langle R_0, R_1\rangle$, with incidence given by non-empty intersection. A map $\cal M$ is compact if and only if $\Omega$ is finite, or equivalently $M$ has finite index in $\Gamma$; it is without boundary and free edges if and only if $M$ is torsion-free (i.e.~contains no conjugate of any $R_i$ or of $R_0R_2$), in which case $\cal M$ is orientable if and only if $M$ is contained in the even subgroup $\Gamma^+$ of $\Gamma$, the subgroup of index $2$ consisting of all words of even length in the generators $R_i$.

The automorphism group $A={\rm Aut}\,{\cal M}$ of $\cal M$ is the set of permutations of $\Omega$ commuting with each $r_i$; equivalently, $A$ is the centraliser of $G$ in ${\rm Sym}\,\Omega$. It acts freely on $\Omega$, since an automorphism fixing a blade must fix all of them, and it is isomorphic to $N_{\Gamma}(M)/M$. We say that $\cal M$ is {\sl regular} if $A$ acts transitively (and hence regularly) on $\Omega$, or equivalently $M$ is normal in $\Gamma$, in which case
\[A\cong G\cong\Gamma/M,\]
and we can regard $A$ and $G$ as the left and right regular representations of the same group. Coverings ${\cal M}_1\to{\cal M}_2$ of maps correspond to inclusions $M_1\leq M_2$ of the corresponding map subgroups, with the index $|M_2:M_1|$ equal to the number of sheets. Normal inclusions correspond to regular coverings, induced by $M_2/M_1$ acting as a group of automorphism of ${\cal M}_2$.

The Petrie polygons of a map $\cal M$ are its closed zig-zag paths, turning alternately first right and first left at successive vertices; these can be identified with the orbits of the dihedral subgroup $\langle R_0R_2, R_1\rangle$ on $\Omega$. Following Coxeter and Moser~\cite[Ch.~8]{CM}, we say that a map $\cal M$ has {\sl type\/} $\{p,q\}_r$, or simply $\{p,q\}$, if $p$, $q$ and $r$ are the orders of the permutations $r_0r_1$, $r_1r_2$ and $r_0r_1r_2$, or equivalently the least common multiples of the valencies of the faces, the vertices and the Petrie polygons of $\cal M$. In a regular map these all have the same valencies $p$, $q$ and $r$ respectively.

In~\cite{Wil}, Wilson described a group $\Sigma\cong S_3$ of operations on regular maps. It is generated by the classical duality $D$, which transposes the vertices and faces of a map, while preserving the underlying surface, and the Petrie duality $P$, which preserves the embedded graph while transposing the faces and the Petrie polygons, so that the surface may change both its orientability and its characteristic. These operations have order $2$. Apart from the identity operation, the remaining elements of $\Sigma$ are the `opposite' operation $DPD=PDP$, also of order $2$, which transposes vertices and Petrie polygons by rejoining adjacent faces across common edges with reversed identifications, and two operations $DP$ and $PD$ of order $3$ called trialities.

In~\cite{JT} it is shown that $\Sigma$ can be identified with the outer automorphism group
\[{\rm Out}\,\Gamma={\rm Aut}\,\Gamma/{\rm Inn}\,\Gamma\cong S_3\]
of $\Gamma$, acting on the normal subgroups $M$ of $\Gamma$ corresponding to regular maps $\cal M$. This is because $\Gamma$ is the free product $V*C$ of a Klein four-group $V=\langle R_0, R_2\rangle\cong V_4$ and a cyclic group $C=\langle R_1\rangle\cong C_2$ of order $2$, and ${\rm Aut}\,\Gamma$ is a semidirect product of the inner automorphism group ${\rm Inn}\,\Gamma\cong\Gamma$ by ${\rm Aut}\,V\cong S_3$, the latter permuting the three involutions $R_0, R_2$ and $R_0R_2$ of $V$. Since the faces, vertices and Petrie polygons of a map $\cal M$ can be identified with the orbits of the dihedral subgroups $\langle R_0, R_1\rangle$, $\langle R_2, R_1\rangle$ and  $\langle R_0R_2, R_1\rangle$ of $\Gamma$ on $\Omega$, one can identify the action of $\Sigma$ on maps with that of ${\rm Aut}\,V$, and hence of ${\rm Out}\,\Gamma$, on normal subgroups $M$ of $\Gamma$.

The direct derivates of a regular map $\cal M$ are its images under $\Sigma$. The only transitive actions of $\Sigma$ are its regular action of degree $6$, its natural action of degree $3$, its action of degree $2$ as the quotient $S_3/A_3\cong S_2$, and its trivial action of degree $1$, so Wilson divided regular maps into classes~I -- IV respectively, as they have $6$, $3$, $2$ or $1$ derivates, admitting one of these four actions. The aim of this paper is to give constructions for the apparently rare maps in class~III: these are the regular maps with just two direct derivates, so that they are invariant under the subgroup $\Sigma^+\cong C_3$ of $\Sigma$ generated by the triality $DP$ (or by its inverse $PD$), but not under $\Sigma$.

If $M$ is a normal subgroup of $\Gamma$, corresponding to a regular map $\cal M$, then let
\[M_{\Sigma}=\bigcap_{\sigma\in\Sigma}\,\sigma(M) \quad{\rm and}\quad
M_{\Sigma^+}=\bigcap_{\sigma\in\Sigma^+}\,\sigma(M),\] 
the largest $\Sigma$-invariant (i.e.~characteristic) and $\Sigma^+$-invariant normal subgroups of $\Gamma$ contained in $M$. Let ${\cal M}_{\Sigma}$ and ${\cal M}_{\Sigma^+}$ be the corresponding maps, the smallest $\Sigma$- and $\Sigma^+$-invariant regular maps covering $\cal M$. Then $\cal M$ is in class~IV if and only if ${\cal M}={\cal M}_{\Sigma}$. or equivalently $M=M_{\Sigma}$, and $\cal M$ is in class~III if and only if ${\cal M}={\cal M}_{\Sigma^+}\neq{\cal M}_{\Sigma}$, or equivalently $M=M_{\Sigma^+}>M_{\Sigma}$.

\medskip

In order to consider the orientability of class~III maps, we need to understand the action of $\Sigma$ on the subgroups $\Delta$ of index $2$ in $\Gamma$. There are seven of these, each containing the commutator subgroup $\Gamma'$, which has quotient $\Gamma/\Gamma'\cong V_8$. Each such subgroup $\Delta$ is determined by which proper subset of the generators $R_0, R_1$ and $R_2$ it contains: we write $\Delta=\Gamma_{\emptyset}$, $\Gamma_i$ or $\Gamma_{ij}$ if it contains none of them, just $R_i$, or just $R_i$ and $R_j$ with $i<j$. Thus $\Gamma_{\emptyset}$ is the even subgroup $\Gamma^+$, consisting of those elements of $\Gamma$ which are words of even length in the generators $R_i$. Both $\Sigma$ and $\Sigma^+$ leave $\Gamma_{02}$ invariant, and have orbits $\{ \Gamma_{\emptyset},\Gamma_0, \Gamma_2\}$ and $\{\Gamma_1, \Gamma_{01},\Gamma_{12}\}$ on the remaining six subgroups. The groups $\Delta$ in the first of these two orbits satisfy $\Delta_{\Sigma}=\Delta_{\Sigma^+}=\Gamma'$, and those in the second satisfy $\Delta_{\Sigma}=\Delta_{\Sigma^+}=\Gamma^*:=\langle\Gamma', R_1\rangle$, subgroups of index $8$ and $4$ in $\Gamma$ respectively.

\begin{lemma}
Each $\Sigma^+$-invariant regular map $\cal M$ has type $\{n,n\}_n$ for some $n$. If $\cal M$ is also orientable and without boundary, then $n$ is even
\end{lemma}

\noindent{\sl Proof.} If $\cal M$ is $\Sigma^+$-invariant, then since $\Sigma^+$ transforms the vertices of $\cal M$ into its faces and Petrie polygons, these all have the same valency $n$, and thus $\cal M$ has type $\{n,n\}_n$. If $\cal M$ is also orientable and without boundary, then the corresponding map subgroup $M$ is contained in $\Gamma^+$. Being $\Sigma^+$-invariant, $M$ is contained $\Gamma^+_{\Sigma^+}=\Gamma'$. Now $R_0R_1$, $R_1R_2$ and $R_0R_1R_2$ have images of order $2$ in $\Gamma/\Gamma'$, and hence the order $n$ of their images in $\Gamma/M\cong {\rm Aut}\,{\cal M}$ is even. \hfill$\square$

\medskip

If a regular map $\cal M$ has type $\{n,n\}_n$ and automorphism group $G$, then it has $|G|/2n$ vertices, $|G|/4$ edges and $|G|/2n$ faces, so it has Euler characteristic
\begin{equation}
\chi=|G|\Bigl(\frac{1}{2n}-\frac{1}{4}+\frac{1}{2n}\Bigr)=|G|\frac{4-n}{4n}.
\end{equation}
In particular, this applies if $\cal M$ is regular and $\Sigma^+$-invariant, in which case all direct derivates of $\cal M$ are on homeomorphic surfaces. These have genus $1-\frac{1}{2}\chi$ or $2-\chi$ as they are orientable or non-orientable.

\section{Generalising Wilson's example}

Wilson's example in~\cite{Wil} is a dual pair $\cal W$ and $D({\cal W})$ of maps of type $\{9,9\}_9$ in class~III with automorphism group $L_2(8)$. Here we will show that there are similar examples based on the group $L_2(2^e)$ whenever $e$ is divisible by $3$. Our constructions use properties of finite fields; see~\cite{LN} for background. 

Let $K$ be the field $F_q$ of order  $q=2^e$, and let $F$ be the subfield $F_r$ where $r=2^f$ and $e=3f$. Then ${\rm Gal}\,K/F$ is a cyclic group of order $3$, generated by the automorphism $\alpha: x\mapsto x^r$. Let us call an element $x\in K$ a {\sl useful generator\/} if it generates the field $K$ and satisfies $x+x^r+x^{r^2}=0$

\begin{lemma}
There exists a useful generator $x\in K$.
\end{lemma}

\noindent{\sl Proof.} The trace map ${\rm Tr}_{K/F}: x\mapsto x+x^r+x^{r^2}$ is a homomorphism $K\to F$ of additive groups, so its kernel has order at least $|K|/|F|=2^{2e/3}$. (In fact we have equality here, since the polynomial ${\rm Tr}_{K/F}(x)$ has at most $r^2=2^{2e/3}$ roots in $K$, but we do not need this.) The elements which do not generate $K$ as a field are those contained in maximal subfields of $K$, and there is one of these, of order $2^{e/p}$, for each prime $p$ dividing $e$. If there are $k$ such primes (including $3$) then $e\geq 3.2^{k-1}$ and so $k\leq 1+\log_2(e/3)$. Each subfield has order at most $2^{e/2}$, so the number of such non-generators is at most $2^{e/2}(1+\log_2(e/3))$. The result is therefore proved if $2^{2e/3}>2^{e/2}(1+\log_2(e/3))$, that is, $2^{e/6}>1+\log_2(e/3)$, and this is true for all $e>6$. The cases $e=3$ and $6$ are easily dealt with by more careful counting arguments (alternatively, see Example~3A). \hfill$\square$

\medskip

The above method of proof has the disadvantages of not being constructive, and of not dealing with small values of $e$. The following constructive proof deals with all cases where $e$ is divisible by $3$ but not by $9$.

\medskip

\noindent{\bf Example 3A.} Suppose that $e=3f$ where $f$ is not divisible by $3$. There exists $t\in E=F_8\leq K$ such that $t^3+t+1=0$. Now $t^7=1$, so $t^{2^i}=t, t^2$ or $t^4$ as $i\equiv 0, 1$ or $2$ mod~$(3)$. Since $r=2^f$ with $f\not\equiv 0$ mod~$(3)$ we have ${\rm Tr}_{K/F}(t)=t+t^r+t^{r^2}=t+t^2+t^4=0$. Now let $x=tu$ where $u$ generates the field $F$. Then ${\rm Tr}_{K/F}(x)=x+x^r+x^{r^2}=tu+(tu)^r+(tu)^{r^2}=(t+t^r+t^{r^2})u=0$ since $u^{r^2}=u^r=u$. The multiplicative order of $u$ divides $2^f-1$, and is therefore coprime to $7$, whereas $t$ has order $7$, so $t$ and $u$ are powers of $x$. Since $t$ and $u$ generate the fields $E$ and $F$, and these subfields generate $K$, it follows that $x$ generates $K$. 

\begin{thm}
For each positive integer $e$ divisible by $3$ there is a regular map in class~III with automorphism group $L_2(2^e)$.
\end{thm}

\noindent{\sl Proof.} Let $q=2^e$ and let $x$ be a useful generator of $K=F_q$. We map the generators $R_i$ of $\Gamma$ into $G=L_2(q)=SL_2(K)$ by

\[R_0\mapsto r_0=\Big(\,\begin{matrix}1&x\cr 0&1\end{matrix}\,\Big),
\quad
R_1\mapsto r_1=\Big(\,\begin{matrix}1&0\cr 1&1\end{matrix}\,\Big)
\quad{\rm and}\quad
R_2\mapsto r_2=\Big(\,\begin{matrix}1&x^r\cr 0&1\end{matrix}\,\Big).
\]
Since $r_0$, $r_1$, $r_2$ and
\[r_0r_2=\Big(\,\begin{matrix}1&x+x^r\cr 0&1\end{matrix}\,\Big)
=\Big(\,\begin{matrix}1&x^{r^2}\cr 0&1\end{matrix}\,\Big)
\]
all have order $2$, this mapping extends to a homomorphism $\theta:\Gamma\to G$. In order to show that $\theta$ is an epimorphism we need to show that no maximal subgroup of $G$ contains $r_0, r_1$ and $r_2$. In~\cite[Ch.~XII]{Dic}, Dickson classified the subgroups of the groups $L_2(q)$ for all prime powers $q$, and in our case we see that the maximal subgroups of $G$ are the following:

\begin{enumerate}

\item normalisers of Sylow $2$-subgroups, isomorphic to the general affine group $AGL_1(q)$;

\item dihedral groups of order $2(q\pm 1)$;

\item subgroups conjugate to $SL_2(q')$ where $q'=2^{e/p}$ for primes $p$ dividing $e$.

\end{enumerate}

\noindent Maximal subgroups of the first type are eliminated because their elements of order $2$ all commute, whereas $r_0$ and $r_1$ do not. A dihedral group of order $2n$, for odd $n$, does not contain a Klein four-group such as $\langle r_0, r_2\rangle$, so maximal subgroups of the second type are also eliminated. Since
\[r_0r_1=\Big(\,\begin{matrix}x+1&x\cr 1&1\end{matrix}\,\Big),\]
has trace $x$ generating $K$, groups of the third type are eliminated. Thus $r_0, r_1$ and $r_2$ generate $G$, so $\theta$ is an epimorphism. The subgroup $M={\rm ker}\,\theta$ of $\Gamma$ is therefore the map subgroup corresponding to a regular map ${\cal M}={\cal M}(x)$ with ${\rm Aut}\,{\cal M}\cong G$.

The automorphism $\alpha$ of order $3$ of $K$ induces an automorphism of $G$ fixing $r_1$ and permuting $r_0$, $r_2$ and $r_0r_2$ in a $3$-cycle $(r_0, r_2, r_0r_2)$, so $\cal M$ is invariant under $\Sigma^+$. If $\cal M$ were invariant under $\Sigma$ then there would be an automorphism of $G$ fixing $r_1$ and transposing $r_0$ and $r_2$. Now ${\rm Aut}\,G\cong \Sigma L_2(K)$, an extension of $G$ by ${\rm Gal}\,K\cong C_e$, and the traces $x$, $x^r$ and $x^{r^2}$ of $r_0r_1$, $r_2r_1$ and $r_0r_2r_1$ each generate the field $K$, so no automorphism of $G$ can transpose $r_0r_1$ and $r_2r_1$ while fixing $r_0r_2r_1$. It follows that $\cal M$ is not invariant under $\Sigma$, so it is in class~III. \hfill$\square$

\medskip

The map ${\cal M}(x)$ constructed here has type $\{n,n\}_n$, where $n$ is the order of the element $r_0r_1$ in $G=L_2(K)$. This can be determined from the trace $x$ of this element. If $K$ is any field of characteristic $2$, then each non-identity element of $L_2(K)$ is conjugate to a matrix
\[A=\Big(\,\begin{matrix}0&1\cr 1&x\end{matrix}\,\Big),\]
where $x$ is its trace. If we write
\[A^n=\Big(\,\begin{matrix}a_n&b_n\cr c_n&d_n\end{matrix}\,\Big),\]
then the equation $A^n=A^{n-1}A$ gives
\[a_n=b_{n-1},\quad b_n=a_{n-1}+xb_{n-1},\quad c_n=d_{n-1},\quad d_n=c_{n-1}+xd_{n-1}.\]
These, and the initial condition $A^0=I$, imply that the sequence $(d_n)$ is the solution of the recurrence relation
\begin{equation}
f_n=f_{n-2}+xf_{n-1}
\end{equation}
satisfying the initial conditions $f_0=1$ and $f_1=x$, and that
\[b_n=c_n=d_{n-1}\quad{\rm and}\quad a_n=d_{n-2}.\]
The sequence defined by $t_n:={\rm Tr}\,A^n=a_n+d_n$ is therefore the solution of the recurrence relation~(2) with initial conditions $f_0=0$ and $f_1=x$. This enables the successive powers of $A$, and in particular their traces, to be computed efficiently, so that one can determine the order of $r_0r_1$ and hence the type of $\cal M$.

\medskip

\noindent{\bf Example 3B.} Let $e=3$. Then $F=F_2$ and we can define $K=F_8$ as $F_2[t]/(t^3+t+1)$, so ${\rm Tr}_{K/F}\,(t)=t+t^2+t^4=0$ and the kernel of ${\rm Tr}_{K/F}$ is $\{0, t, t^2, t^4\}$. Taking $x=t$ as a useful generator for $K$, the proof of Theorem~3.2, gives
\[r_0r_1=\Big(\,\begin{matrix}t+1&t\cr 1&1\end{matrix}\,\Big).\]
Applying the recurrence relation~(1) we find that $(r_0r_1)^3$ has trace $t_3=1$, so it has order $3$ and hence $r_0r_1$ has order $9$. Thus the corresponding map ${\cal M}={\cal M}(x)$ and its dual $D({\cal M})$ have type $\{9,9\}_9$, so by equation~(1) they have characteristic $\chi=-70$.

\medskip 

As we will now show, this map $\cal M$ is isomorphic to the class~III map $\cal W$ described by Wilson in~\cite{Wil}. His construction is completely different from ours, and it does not seem to generalise so easily to give other class~III maps. He starts with the non-orientable regular map $\cal N$ denoted in~\cite{CM} by $\{3,7\}_9$, since it is the largest regular map of this type. The opposite map $DPD({\cal N})$ is the non-orientable regular map $\{3,9\}_7$ of this corresponding type. Having defined an operation $H_2$, applicable to maps of odd vertex valency, which preserves the embedded graph but replaces the rotation $r_1r_2$ of edges around each vertex with its square, he notes that the regular map ${\cal W}=H_2(DPD({\cal N}))$ has type $\{9,9\}_9$. He then shows that ${\rm Aut}\,{\cal W}$ has generators $R$, $S$ and $T$ (corresponding to our $r_0r_1$, $r_2r_1$ and $r_0r_2r_1$), with defining relations
\begin{equation}
R^9=S^9=T^9=(RS^3)^3=(ST^3)^3=(TR^3)^3=(SR^3)^7=(TS^3)^7=(RT^3)^7=1.
\end{equation}
The $3$-cycle $(R,S,T)$ clearly induces an automorphism of ${\rm Aut}\,{\cal W}$ whereas the transposition $(R,S)$ does not, so $\cal W$ is in class~III. Wilson does not specifically identify ${\rm Aut}\,{\cal W}$, but since $D$, $P$ and $H_2$ preserve automorphism groups we have ${\rm Aut}\,{\cal W}={\rm Aut}\,{\cal N}$, and the latter is identified by Coxeter and Moser in~\cite[Table~8]{CM} as $LF(2,2^3)$, Dickson's notation for $L_2(8)$. It is straightforward to check that the generators $r_i$ in Example~3B satisfy the relations corresponding to $(3)$, so ${\cal M}\cong{\cal W}$.

Alternatively, one can identify the maps in Example~3B with Wilson's pair by using the character table and list of maximal subgroups of $L_2(8)$ in~\cite{CCNPW}, together with formula~(7.3) of~\cite{Ser}, to show that the triangle group
\[\Delta=\langle X, Y, Z\mid X^9=Y^2=Z^9=XYZ=1\rangle\]
has only two normal subgroups with quotient group $L_2(8)$. These lift back, through the epimorphism $R_1R_2\mapsto X$, $R_2R_0\mapsto Y$, $R_0R_1\mapsto Z$, to two normal subgroups of $\Gamma^+$ (which are, in fact, normal in $\Gamma$) corresponding to the canonical double covers of two regular maps of type $\{9,9\}$ with automorphism group $L_2(8)$. Any such map must arise in this way, so there are only two of them, namely $\cal W$ and its dual.

One can also verify this identification by using Conder's lists of regular maps~\cite{Con}. He used a computer to find all the orientable and non-orientable regular maps of characteristic $\chi=-1, \ldots, -200$ (those with $\chi\geq 0$ are well-known, see~\cite[Ch.~8]{CM} for instance, and there are none in class~III). In his notation, an entry R$g.i$ or N$g.i$ denotes the $i$-th orientable or non-orientable regular map (or dual pair of maps) of genus $g$, lexicographically ordered by their type $\{p,q\}$ with $p\leq q$. There is only one entry corresponding to regular maps of characteristic $-70$ and type $\{9,9\}$, namely the dual pair N72.9 of type $\{9,9\}_9$; their canonical double covers are the pair R71.15 of type $\{9,9\}_{18}$. The map ${\cal N}=\{3,7\}_9$ used by Wilson is Conder's N8.1, while $DPD({\cal N})$ is N16.1.

\begin{thm}
{\rm(a)} The only regular maps in class~III with automorphism group $L_2(q)$ for $q=2^e\geq 4$ are the maps ${\cal M}(x)$ constructed in Theorem~3.2, together with their dual maps $D({\cal M}(x))$, where $e$ is divisible by $3$ and $x$ is a useful generator for the field $K=F_q$.
\vskip2pt
\noindent{\rm(b)} Two such maps ${\cal M}(x_1)$ and ${\cal M}(x_2)$ are isomorphic if and only if $x_1$ and $x_2$ are conjugate under ${\rm Gal}\,K$.
\vskip2pt
\noindent{\rm(c)} There are no  isomorphisms between maps ${\cal M}(x_1)$ and dual maps $D({\cal M}(x_2))$.  
\end{thm}

\noindent{\sl Proof.} (a) Let $\cal M$ be a regular map in class~III with ${\rm Aut}\,{\cal M}\cong G=L_2(q)$, corresponding to an epimorphism $\Gamma\to G$, $R_i\mapsto r_i$.
Since $G$ is neither dihedral nor cyclic, each $r_i$ is an involution. All involutions in $G$ are conjugate, so by composing $\theta$ with an inner automorphism of $G$ we may assume that
\[ r_1=\Big(\,\begin{matrix}1&0\cr 1&1\end{matrix}\,\Big). \]
Now $r_0$ and $r_2$ are commuting involutions. The centraliser of any involution in $G$ is the unique Sylow $2$-subgroup containing it, so $r_0$ and $r_2$ lie in a Sylow $2$-subgroup $S$, which must be distinct from the Sylow $2$-subgroup $T=C_G(r_1)$ containing $r_1$. Each Sylow $2$-subgroup acts regularly by conjugation on the $q$ others in $G$, so by further conjugating with a unique element of $T$ we may also assume that
\[r_0=\Big(\,\begin{matrix}1&x\cr 0&1\end{matrix}\,\Big)
\quad{\rm and}\quad
r_2=\Big(\,\begin{matrix}1&y\cr 0&1\end{matrix}\,\Big),
\quad{\rm giving}\quad
r_0r_2=\Big(\,\begin{matrix}1&x+y\cr 0&1\end{matrix}\,\Big),\]
for some $x, y\in K$. Indeed, this argument shows that there is a {\sl unique\/} inner automorphism of $G$ taking the original generators $r_0, r_1$ and $r_2$ to this form.

Since $\cal M$ is in class~III there is an automorphism $\alpha$ of $G$ fixing $r_1$ and permuting $r_0, r_2$ and $r_0r_2$ in a $3$-cycle. Now ${\rm Aut}\,G=\Sigma L_2(q)$, a semidirect product of $G$ by ${\rm Gal}\,K$. The automorphisms fixing $r_1$ form the semidirect product of $T$ by ${\rm Gal}\,K$, and the above uniqueness result shows that such an automorphism preserves $S$ if and only if it is in ${\rm Gal}\,K$. Since $\alpha$ has order $3$ it follows that $e=3f$ for some $f$, and $\alpha$ acts on $K$ by $t\mapsto t^r$ or $t\mapsto t^{r^2}$, where $r=2^f$. Suppose first that $\alpha: t\mapsto t^r$. Then $y=x^r$ and $x+y=x^{r^2}$, so $x+x^r+x^{r^2}=0$. If $K_0$ is the subfield of $K$ generated by $x$ then $y\in\alpha(K_0)=K_0$, so $r_0, r_1, r_2\in L_2(K_0)$; however, these elements generate $G$, so $K_0=K$. Thus $x$ is a useful generator of $K$, and ${\cal M}\cong{\cal M}(x)$. We obtain a similar conclusion if $\alpha:t\mapsto t^{r^2}$, except that now $y=x^{r^2}$ and $x+y=x^r$; this differs from the preceding case by inverting $\alpha$, or equivalently by transposing $r_0$ and $r_2$, so in this case ${\cal M}\cong D({\cal M}(x))$.

\smallskip

\noindent(b) This condition is sufficient, since an automorphism of $K$ taking $x_1$ to $x_2$ will give an automorphism of $G$ inducing an isomorphism from ${\cal M}(x_1)$ to ${\cal M}(x_2)$. The uniqueness result in the proof of (a) also shows that this condition is necessary.

\noindent(c) The set of maps ${\cal M}(x)$ is distinguished from the set of their dual maps by the fact that the field automorphism inducing the triality sending $r_0$ to $r_2$ is respectively given by $t\mapsto t^r$ or its inverse, so there can be no isomorphisms between these two sets of maps. \hfill$\square$

\medskip

Since ${\rm Gal}\,K$ has order $e$, and acts wthout fixed points of the set of useful generators of $K$, we immediately have the following corollary:

\begin{cor}
If $q=2^e$ with $e$ divisible by $3$, and $N_e$ is the number of useful generators of $K=F_q$, then there are, up to isomorphism, $2N_e/e$ regular maps in class~III with automorphism group $L_2(q)$, forming $N_e/e$ dual pairs. \hfill$\square$
\end{cor}

\noindent{\bf Example 3C.} Let $e=6$, so that $G=L_2(64)$. The field $K=F_{64}$ has two maximal subfields $E=F_8=F_2[t]/(t^3+t+1)$ and $F=F_4=F_2[u]/(u^2+u+1)$. There are $2^{2e/3}=16$ elements $x$ in the kernel $Z$ of the homomorphism ${\rm Tr}_{K/F}:x\mapsto x+x^4+x^{16}$. If $x\in F$ then $x=x^4=x^{16}$, so only $x=0$ is in $Z$. If $x\in E$ then $x^8=x$ so $x^{16}=x^2$ and hence ${\rm Tr}_{K/F}(x)=x+x^2+x^4={\rm Tr}_{E/F_2}(x)$, giving $2^2=4$ elements of $Z$ (namely $0, t, t^2$ and $t+t^2$). Thus there are $N_6=16-4=12$ useful generators $x$ for $K$, giving $N_6/6=2$ dual pairs of class~III maps with automorphism group $G$.

In order to construct such a map we can take $x=tu$ as a useful generator for $K$, as in Example~3A. Applying the recurrence relation~(2), and using $t^3=t+1$ and $u^2=u+1$, we eventually find that $(r_0r_1)^{13}$ has trace $t_{13}=u+1$, so it is conjugate to an element of order $5$ in $L_2(F)$; since $t_5=t^2+u+t^2u\neq 0$ we have $(r_0r_1)^5\neq I$, so it follows that $r_0r_1$ has order $65$. Thus the resulting dual pair of class~III maps have type $\{65,65\}_{65}$, and since they have automorphism group $G$ of order $63.64.65$, they have characteristic $\chi=-61488$ by equation~(1). Alternatively, since ${\rm Tr}_{K/F}(t^2)={\rm Tr}_{K/F}(tu)=0$ we can take $x=t^2+tu$ as a useful generator. In this case $t_7=t^2+t$, which is the trace of an element of order $9$ in $L_2(E)$, so $r_0r_1$ has order $63$. This gives the second dual pair of class~III maps, of type $\{63,63\}_{63}$ and characteristic $-61360$.

\medskip

We can generalise the counting argument in Example~3C, as follows. Let $Z$ be the kernel of ${\rm Tr}_{K/F}$, the set of zeros of $x+x^r+x^{r^2}$, so $|Z|=2^{2e/3}$. An element $x\in K$ is a useful generator for $K$ if it is in $Z$ but not in any proper subfield $L$ of $K$. The automorphism $\alpha:x\mapsto x^r$ of $K$ has fixed field $F$, so ${\rm Tr}_{K/F}(x)=x$ for all $x\in F$; thus if $L$ is a subfield of $F$ then $L\cap Z=\{0\}$. If $L$ is not a subfield of $F$, then $\alpha$ induces an automorphism of order $3$ of $L$ with fixed field $M=L\cap F$, so ${\rm Tr}_{K/F}(x) = {\rm Tr}_{L/M}(x)$ for all $x\in L$,  giving $|L\cap Z|=|L|/|M|=|L|^{2/3}$. Now $K$ has one subfield $L$, of order $2^d$, for each $d$ dividing $e$, and $L$ is a subfield of $F$ if and only if $d$ divides $f$. M\"obius inversion over the lattice of subfields of $K$, partitioned into those contained or not contained in $F$, gives
\[N_e=\sum_{d\mid e}\mu\Bigl(\frac{e}{d}\Bigr)|L\cap Z|=
\sum_{d\mid e,\, d\negthinspace\not\;\mid f}\mu\Bigl(\frac{e}{d}\Bigr)2^{2d/3}
+\sum_{d\mid f}\mu\Bigl(\frac{e}{d}\Bigr)\]
where $\mu$ is the M\"obius function~\cite[\S\S 16.3--4]{HW}. Putting $c=e/d$ we can write this as
\[N_e=\sum_{c\mid e'}\mu(c)2^{2e/3c}+\sum_{3\mid c\mid e}\mu(c)\]
where $e=3^ie'$ with $(e',3)=1$. Putting $c=3b$ in the second sum gives
\[\sum_{3\mid c\mid e}\mu(c)=\sum_{b\mid f}\mu(3b)
=-\sum_{3\negthinspace\not\;\mid b\mid f}\mu(b)
=-\sum_{b\mid e'}\mu(b),\]
and this is $-1$ or $0$ as $e'=1$ or $e'>1$, that is, as $e$ is or is not a power of $3$. We have therefore proved:

\begin{prop}
Let $e=3^ie'$ with $e'$ coprime to $3$. Then
\[N_e=\sum_{c\mid e'}\mu(c)2^{2e/3c}-\nu_e\]
where $\mu$ is the M\"obius function and $\nu_e=1$ or $0$ as $e'=1$ or $e'>1$.
\hfill$\square$
\end{prop}

\noindent{\bf Example 3D.} Together with Corollary~3.4, this shows that if $e=3, 6, 9, 12$ or $15$ then the number of dual pairs of class~III maps with automorphism group $L_2(2^e)$ is respectively $1, 2, 7, 20$ or $68$. If $e=3^i>3$ we have $N_e=2^{2e/3}-1$, giving $(2^{2e/3}-1)/e$ such pairs.

\medskip

The argument used to prove Lemma~3.1 shows that the formula in Proposition~3.5 is dominated by the term with $c=1$. Thus $N_e\sim 2^{2e/3}=q^{2/3}$ as $e\to\infty$, so the number $N_e/e$ of dual pairs in Corollary~3.4 is asymptotic to $q^{2/3}/e$, which grows quite rapidly. This shows that one can find arbitrarily many class~III pairs with the same automorphism group $G=L_2(2^e)$ by taking $e$ sufficiently large. These maps have type $\{p,p\}_p$ where $p$ is the order of $r_0r_1$; now the orders of the elements of $G$ are $2$ and the divisors of $q\pm 1$, so the number of distinct orders is $d(q-1)+d(q+1)$ where $d(n)$ is the number of divisors of an integer $n$. As $n\to\infty$ this function grows more slowly than any positive power of $n$ (see~\cite[Theorem 315]{HW}, for instance), so by taking $e$ sufficiently large one can find arbitrarily many class~III pairs with automorphism group $G$, and with mutually distinct types.


\section{Maps in class III as coverings}

Another method of finding maps in class~III is to construct them as coverings of known maps in class~III: for instance, several of the examples found by Conder in~\cite{Con06} are coverings of Wilson's maps $\cal W$ and $D({\cal W})$. We start with a normal subgroup $N$ of $\Gamma$, corresponding to a regular map $\cal N$ in class~III, and we try to find subgroups $M$ of $N$, corresponding to maps $\cal M$ which cover $\cal N$, such that $\cal M$ is also in class~III.

\begin{lemma}
Let $N$ be the map subgroup of $\Gamma$ corresponding to a map in class~III, and let $M$ be a characteristic subgroup of $N$ such that $N/M$ is a characteristic subgroup of $\Gamma/M$. Then the map $\cal M$ corresponding to $M$ is in class~III.
\end{lemma} 

\noindent{\sl Proof.} Since $M$ is a characteristic subgroup of $N$, and $N$ is a $\Sigma^+$-invariant normal subgroup of $\Gamma$, it follows that $M$ is a $\Sigma^+$-invariant normal subgroup of $\Gamma$. The corresponding map $\cal M$ is therefore regular and invariant under $\Sigma^+$. If $\Sigma$ leaves $M$ invariant, it induces a group of automorphisms of $\Gamma/M$ which must leave the characteristic subgroup $N/M$ invariant. Thus $\Sigma$ must leave $N$ invariant, which is false, so $M$ is not $\Sigma$-invariant and the corresponding map $\cal M$ is in class~III. \hfill$\square$

\medskip

Now $N$ is the fundamental group of a punctured surface, namely the underlying surface of $\cal N$, with punctures at the vertices and face centres, so it is a free group of finite rank $r$ (equal to $V+F+2g-1$ if $\cal N$ is orientable, of genus $g$, and $V+F+g-1$ if it is non-orientable of genus $g$). For any integer $n\geq 2$, the subgroup $M=N'N^n$ of $N$ generated by its commutators and $n$th powers is a characteristic subgroup of index $n^r$ in $N$. In most cases, and certainly if $n$ is coprime to $|\Gamma:N|$, $N/M$ is a characteristic subgroup of $\Gamma/M$, so Lemma~4.1 gives a class~III map $\cal M$ corresponding to $M$. In cases where ${\rm Aut}\,{\cal N}$ is a non-abelian simple group (as when $\cal N$ is Wilson's map $\cal W$, for instance), this is true for all $n$, since $N/M$ is the only normal subgroup of $\Gamma/M$ with a non-abelian simple quotient group.

If we replace $\Gamma$ with the appropriate extended triangle group by adding the relations $(R_0R_1)^p=(R_1R_2)^q=1$, where $\cal N$ has type $\{p,q\}$, then $N$ becomes a surface group, of rank $r=2g$ or $g$ as $\cal N$ is orientable or not, and a similar construction can be used.

Choosing different values of $n$ gives infinitely many class~III maps $\cal M$ as coverings of $\cal N$. Alternatively one could iterate this process, finding a suitable characteristic subgroup of $M$, and so on. Taking $M$ to be a subgroup of infinite index, such as $N'$, gives infinite maps in class~III by the same argument.

Another way of constructing class~III maps as coverings is to start with the map subgroup $N$ corresponding to a class~III map $\cal N$, and to take its intersection $M=N\cap K$ with a $\Sigma^+$-invariant normal subgroup $K$ of $\Gamma$ such that $N/M$ is a characteristic subgroup of $\Gamma/M$. Since $N$ and $K$ are normal and $\Sigma^+$-invariant, so is $M$. If $M$ were $\Sigma$-invariant then $\Sigma$ would induce a group of automorphisms of $\Gamma/M$ preserving $N/M$, so $N$ would be $\Sigma$-invariant, which is false; thus $M$ is not $\Sigma$-invariant, so it corresponds to a map $\cal M$ in class~III.  For example, if $K=\Gamma_{02}$, $\Gamma^*$ or $\Gamma'$, of index $2$, $4$ or $8$, so that $N/M\cong NK/K$ is an elementary abelian $2$-group, and if ${\rm Aut}\,{\cal N}$ has no non-trivial normal $2$-subgroups, then $N/M$ is a characteristic subgroup of $\Gamma/M$, as required.

\medskip

\noindent{\bf Example 4A.} By taking ${\cal N}={\cal W}$ or $D({\cal W})$ and $K=\Gamma_{02}$ we obtain a pair of maps $\cal M$ in class~III which are non-orientable double covers of Wilson's maps; these are the dual pair N198.6 of type $\{18,18\}_{18}$ and characteristic $-196$ in~\cite{Con}, with automorphism group $L_2(8)\times C_2$. If, instead, we take $K=\Gamma^*$ we obtain a dual pair of non-orientable maps of type $\{18,18\}_{18}$ and characteristic $-392$ with automorphism group $L_2(8)\times V_4$, whereas taking $K=\Gamma'$ gives a dual pair of orientable maps of type $\{18,18\}_{18}$ and characteristic $-784$ with automorphism group $L_2(8)\times V_8$.

\medskip

On the other hand, $\Gamma^+$ is not $\Sigma^+$-invariant, so the canonical double cover of a non-orientable class~III map, corresponding to $M=N\cap\Gamma^+$, need not be in class~III: Again, Wilson's maps $\cal W$ and $D({\cal W})$ illustrate this, having canonical double covers in class~I, namely the dual pair R71.15 in~\cite{Con}.

\section{Maps in class III as parallel products}

Let $\cal N$ be a regular map in class~I, corresponding to a normal subgroup $N$ of $\Gamma$. Since $\cal N$ has six direct derivates, $N$ is one of an orbit $\{N_1,\ldots, N_6\}$ of six normal subgroups of $\Gamma$ under the action of $\Sigma$, all with quotient group $\Gamma/N_i\cong G={\rm Aut}\,{\cal N}$. We can choose the numbering so that $N_1=N$, and $N_2$ and $N_3$ are the remaining images of $N$ under $\Sigma^+$. The subgroup $M=N_{\Sigma^+}=N_1\cap N_2\cap N_3$ is then normal in $\Gamma$, and has finite index in $\Gamma$ provided $N$ has finite index, so it corresponds to a regular map ${\cal M}={\cal N}_{\Sigma^+}$ with ${\rm Aut}\,{\cal M}\cong\Gamma/M$. In Wilson's terminology~\cite{Wil94}, $\cal M$ is the parallel product of the maps $\cal N$, $DP({\cal N})$ and $PD({\cal N})$; if $\cal N$ has type $\{p,q\}_r$ then $\cal M$ has type $\{n,n\}_n$ where $n={\rm lcm}\{p,q,r\}$. Clearly $M$ is $\Sigma^+$-invariant, and it is $\Sigma$-invariant if and only if it coincides with $N_{\Sigma}=N_1\cap\cdots\cap N_6$, so $\cal M$ is in class~III or IV as $M>N_{\Sigma}$ or $M=N_{\Sigma}$ respectively. There are many cases where the former happens, and the following simple (and presumably well-known) lemma provides a sufficient condition for this.

\begin{lemma} Let $H$ be a group with distinct normal subgroups $K_1,\ldots, K_k$ such that each $H/K_i$ is isomorphic to a non-abelian simple group $G_i$. Then the normal subgroup $K=K_1\cap\cdots\cap K_k$ of $H$ has quotient $H/K\cong G_1\times\cdots\times G_k$.
\end{lemma}

\noindent{\sl Proof.} We use induction on $k$. There is nothing to prove if $k=1$, so assume that $k>1$ and that the result has been proved for $k-1$ normal subgroups. Let $L=K_1\cap\cdots\cap K_{k-1}$. We first need to show that $L\not\leq K_k$, so suppose that $L\leq K_k$. Then $K_k/L$ is a maximal normal subgroup of $H/L$; however, $H/L\cong G_1\times\cdots\times G_{k-1}$ by the induction hypothesis, and since the groups $G_i$ are non-abelian and simple the maximal normal subgroups of $G_1\times\cdots\times G_{k-1}$ are just the $k-1$ products of all but one of the direct factors~\cite[I.9.12(b)]{Hup}, corresponding to the subgroups $K_i/L$ of $H/L$ for $i=1,\ldots, k-1$; thus $K_k/L=K_i/L$ and hence $K_k=K_i$ for some $i<k$, against the hypotheses. Thus $L\not\leq K_k$, so $H=K_kL$ since $K_k$ is a maximal normal subgroup of $H$; since $K=L\cap K_k$ we therefore have
\[H/K = K_k/K\times L/K \cong H/L\times H/K_k \cong 
(G_1\times\cdots\times G_{k-1})\times G_k,\]
as required. \hfill$\square$

\medskip

This result cannot be extended to arbitrary quotient groups $G_i$: for instance, if $p$ is prime and $e>1$ then an elementary abelian group $H=C_p\times\cdots\times C_p$ of order $p^e$ has $k=(p^e-1)/(p-1) > e$ normal subgroups $K_i$ with quotient $H/K_i\cong C_p$ and intersection $K=1$. In general, all one can conclude is that the projections $H/K\to H/K_i$ induce an embedding of $H/K$ in $G_1\times\cdots\times G_k$.

A group $A$ is {\sl almost simple\/} if it has a non-abelian simple normal subgroup $S$ such that $C_A(S)=1$, so that $A$ is embedded in ${\rm Aut}\,S$. Then $S$ is the unique minimal normal subgroup of $A$, since any other would centralise $S$. If $A$ is finite then $A/S$ is solvable by the Schreier Conjecture, which asserts that if $S$ is any non-abelian finite simple group then ${\rm Out}\,S$ is solvable; the general proof depends on the classification of finite simple groups (see~\cite[Theorem~1.46]{Gor}, for instance), though in any specific case it can be verified directly.

\begin{thm}
If $\cal N$ is a regular map in class~I such that ${\rm Aut}\,{\cal N}$ is almost simple then the map ${\cal M}={\cal N}_{\Sigma^+}$ is in class~III.
\end{thm}

\noindent{\sl Proof.} Let $\cal N$ correspond to an epimorphism $\theta:\Gamma\to A={\rm Aut}\,{\cal N}$ with kernel $N$, and let $M=N_{\Sigma^+}$. It is sufficient to show that $M\neq N_{\Sigma}$, so suppose for a contradiction that $M=N_{\Sigma}$. Let $T=\theta^{-1}(S)$, where $S$ is the unique non-abelian simple normal subgroup of $A$, and let $N_i$ and $T_i$, for $i=1,\ldots, 6$, be the images of $N$ and $T$ under $\Sigma$, numbered as above. Thus $T_i$ is a normal subgroup of $\Gamma$ containing $N_i$, with $\Gamma/T_i\cong A/S$ solvable and $T_i/N_i\cong S$ non-abelian and simple.  Now $T_{\Sigma^+}$ and $T_{\Sigma}$ are the smallest normal subgroups of $\Gamma$, respectively containing $N_{\Sigma^+}$ and $N_{\Sigma}$, with solvable quotient groups, so they are equal since $N_{\Sigma^+}=N_{\Sigma}$. Let $H=T_{\Sigma^+}=T_{\Sigma}$. The six normal subgroups $K_i=H\cap N_i$ of $H$ are all distinct, since if $K_i=K_j$ with $i\neq j$ then $N_iN_j/N_i$ is a non-trivial solvable normal subgroup of the almost simple group $\Gamma/N_i$, which is impossible. Now $M$ is the intersection of both three and six of these groups $K_i$, and they have non-abelian simple quotient groups $H/K_i\cong HN_i/N_i=T_i/N_i\cong S$, so Lemma~5.1 gives $S^3\cong H/M\cong S^6$, which is impossible. Thus $M\neq N_{\Sigma}$, so the corresponding map ${\cal M}={\cal N}_{\Sigma^+}$ is in class~III. \hfill$\square$

\medskip

If $\cal N$ is orientable then $\cal M$, which covers $\cal N$, is also orientable. If $A=S$, that is, ${\rm Aut}\,{\cal N}$ is simple, then $H=\Gamma$ and the above proof shows that ${\rm Aut}\,{\cal M} = S^3$. However, if $A>S$ then all one can conclude is that $S^3\unlhd {\rm Aut}\,{\cal M}\leq A^3$, and in general more work is required to identify ${\rm Aut}\,{\cal M}$ precisely (see Examples~5A and 5C).

\medskip

There are many instances of maps $\cal N$ satisfying the hypotheses of Theorem~5.2, and here we give a few straightforward examples.

\medskip

\noindent{\bf Example 5A.} The smallest orientable map $\cal N$ satisfying the hypotheses of Theorem~5.2 is the map R3.1 in~\cite{Con}, of type $\{3,7\}_8$ and characteristic $-4$; it has automorphism group $A=PGL_2(7)$ with $S=L_2(7)$, the simple group of order $168$. (The underlying surface of $\cal N$ can be identified with Klein's quartic curve of genus $3$, the Riemann surface of least genus attaining Hurwitz's upper bound of $84(g-1)$ automorphisms for surfaces of genus $g\geq 2$.) Theorem~5.2 then gives an orientable class~III map $\cal M$ of type $\{168,168\}_{168}$. Each $r_i$ is in $A\setminus S$, so $T=\theta^{-1}(S)$ is the group $\Gamma^+=\Gamma_{\emptyset}$ and hence $H=T_{\Sigma^+}=\Gamma'$ (see \S 2). The proof of Theorem~5.2 then shows that ${\rm Aut}\,{\cal M}=A^3=PGL_2(7)^3$. Since $|PGL_2(7)|=336$, equation~(1) shows that $\cal M$ has characteristic $-9257472$. For the smallest non-orientable map $\cal N$ satisfying the hypotheses of Theorem~5.2, see the case $n=5$ of Example~5C.


\medskip

\noindent{\bf Example 5B.} Let $A$ be the simple group $L_2(p)=SL_2(p)/\{\pm 1\}$ for a prime $p\equiv 1$ mod~$(4)$. Since $-1=i^2$ for some $i\in F_p$ we can define a homomorphism $\theta:\Gamma\to A$ by
\[R_0\mapsto r_0=\pm\Big(\,\begin{matrix}0&i\cr i&0\end{matrix}\,\Big),
\quad
R_1\mapsto r_1=\pm\Big(\,\begin{matrix}i&i\cr 0&-i\end{matrix}\,\Big)
\quad{\rm and}\quad
R_2\mapsto r_2=\pm\Big(\,\begin{matrix}i&0\cr 0&-i\end{matrix}\,\Big).
\]
Thus $\langle r_0, r_2\rangle$ is a Klein four-group, and
\[
r_1r_2=\pm\Big(\,\begin{matrix}1&-1\cr 0&1\end{matrix}\,\Big),
\]
has order $p$. The maximal subgroups of $A$ are isomorphic to the unique subgroup of order $p(p-1)/2$ in $AGL_1(p)$, to dihedral groups of order $p\pm 1$, and (depending on $p$) to $A_4$, $S_4$ or $A_5$. If we take $p>5$ then only the first of these contain elements of order $p$, and these do not contain any Klein four-groups, so $r_0, r_1$ and $r_2$ generate $A$. Thus $\theta$ is an epimorphism, so its kernel $N$ corresponds to a regular map $\cal N$ with ${\rm Aut}\,{\cal N}\cong A$. Now
\[
r_0r_1=\pm\Big(\,\begin{matrix}0&1\cr -1&-1\end{matrix}\,\Big)
\quad{\rm and}\quad
r_0r_2r_1=\pm\Big(\,\begin{matrix}0&i\cr i&i\end{matrix}\,\Big),
\]
with $i\neq\pm 1$ or $\pm 2$ for $p>5$, so these and $r_1r_2$ all have distinct traces. Since ${\rm Aut}\,A=PGL_2(p)$ preserves traces, it follows that $\cal N$ is in class~I. We therefore obtain a non-orientable regular map $\cal M$ in class~III with ${\rm Aut}\,{\cal M}\cong L_2(p)^3$ for each prime $p>5$ satisfying $p\equiv 1$ mod~$(4)$. Since $r_0r_1$ and $r_1r_2$ have orders $3$ and $p$, $\cal N$ has type $\{3,p\}_r$ where $r$ is the order of $r_0r_1r_2$. It follows that $\cal M$ has type $\{n,n\}_n$ where $n$ is $pr$ or $3pr$ as $3$ does or does not divide $r$. In the smallest example, arising for $p=13$, we find that $r=7$, so $\cal M$ has type $\{273,273\}_{273}$ and characteristic $-320772816$; the map $\cal N$ is N51.1 in~\cite{Con}.

(The canonical double covers of the maps $\cal N$ used here are the regular triangular maps associated with the principal congruence subgroups of level $p$ in the modular group $PSL_2({\bf Z})$; these orientable maps have automorphism group $L_2(p)\times C_2$ for $p\equiv 1$ mod~$(4)$.)

\medskip

\noindent{\bf Example 5C.} Let us define a homomorphism $\theta:\Gamma\to S_n$ for $n\geq 5$ by
\[R_0\mapsto r_0=(1,n)(2,n-1)(3,n-2)\ldots,\quad R_1\mapsto r_1=(2,n)(3,n-1)(4,n-2)\ldots\]
and
\[R_2\mapsto r_2=(1,n).
\]
Here $r_0$ and $r_1$ are induced by two reflections of a regular $n$-gon with vertices $1, 2, \ldots, n$, while $r_2$ is chosen to commute with $r_0$. Then
\[r_0r_1= (1,2,3,\ldots, n),\quad r_1r_2 = (1,n,2)(3,n-1)(4,n-2)\ldots\]
and
\[r_0r_1r_2 = (1, 2,3,\ldots,n-1),\]
these permutations having orders $n$, $6$ and $n-1$. The $n$-cycle $r_0r_1$ and the transposition $r_2$ generate $S_n$, so $\theta$ is an epimorphism. The subgroup $N=\ker\theta$ of $\Gamma$ therefore corresponds to a regular map $\cal N$ of type $\{n,6\}_{n-1}$ with ${\rm Aut}\,{\cal N}\cong S_n$. If $n\neq 6$ then all automorphisms of $S_n$ are inner; now $r_0$, $r_2$ and $r_0r_2$ have different cycle-structures if $n\geq 7$, so they are mutually inequivalent under automorphisms of $S_n$, and hence $\cal N$ is in class~I. The same applies for $n=5$ since then $\cal N$ has type $\{5,6\}_4$. If $n=6$, however, $\cal N$ is the map N62.3 of type $\{6,6\}_5$ in~\cite{Con}, and the outer automorphism of $S_6$, interchanging transpositions with products of three transpositions, induces a self-duality. Since $S_n$ is almost simple for $n\geq 5$ we therefore obtain a regular map ${\cal M}={\cal N}_{\Sigma^+}$ in class~III for $n=5$ and for each $n\geq 7$. 

Note that $r_0$ is even if and only if $n\equiv 0$ or $1$ mod~$(4)$, that $r_1$ is even if and only if $n\equiv 1$ or $2$ mod~$(4)$, and that $r_2$ is always odd. Thus the subgroup $T=\theta^{-1}(A_n)$ is $\Gamma_0$, $\Gamma_{01}$, $\Gamma_1$ or $\Gamma_{\emptyset}$ as $n\equiv 0, 1, 2$ or $3$ mod~$(4)$. It follows from \S 2 that the subgroup $H=T_{\Sigma^+}$ is $\Gamma'$ if $n\equiv 0$ or $3$ mod~$(4)$, and it is $\Gamma^*=\langle\Gamma', R_1\rangle$ if $n\equiv 1$ or $2$ mod~$(4)$. In the first case $\cal M$ is orientable, with ${\rm Aut}\,{\cal M}=S_n^3$, whereas in the second case it is non-orientable, with ${\rm Aut}\,{\cal M}$ a subgroup of index $2$ in $S_n^3$. The smallest example, for $n=5$, is a non-orientable map of type $\{60,60\}_{60}$ and characteristic $-201600$; in this case $\cal N$ is a direct derivate of the map N5.1 of type $\{4,5\}_6$, an antipodal quotient of the map R4.2 of this type on Bring's curve of genus $4$. (In fact $\cal N$ and its direct derivates are the smallest non-orientable maps satisfying the hypotheses of Theorem~5.2.)

\medskip

\noindent{\bf Example 5D.} In Example~5C let $n\equiv 1$ mod~$(4)$, so that $r_0$ and $r_1$ are even. There are several choices of even permutations $r_2$ giving rise to class~I maps $\cal N$ with ${\rm Aut}\,{\cal N}=A_n$, and hence to class~III maps $\cal M$ with ${\rm Aut}\,{\cal M}=A_n^3$. For example, if $n>5$ and
\[r_2=(1,n)(3,n-2)\]
then $r_1r_2=(1,n,2)(3,n-1,n-2,4)(5,n-3)\ldots$, so the elements $(r_1r_2)^4=(1,n,2)$ and $r_0r_1=(1,2,\ldots, n)$ generate $A_n$. Since $r_1r_2$ and $r_0r_1r_2=(1,2,n-2,n-1)(3,4,\ldots, n-3)$ have orders $12$ and $n-5$, $\cal N$ has type $\{n,12\}_{n-5}$. Thus $\cal N$ is in class~I for $n\neq 17$, and if $n=17$ this follows from the fact that $r_0, r_2$ and $r_0r_2$ have different cycle structures. Since $A_n$ is simple for $n\geq 5$, Theorem~5.2 gives a non-orientable class~III map $\cal M$ with ${\rm Aut}\,{\cal M}=A_n^3$. The smallest example, with $n=9$, has type $\{36,36\}_{36}$ and characteristic $-1327353495552000$.

As an alternative, if we take $r_2=(1,2)(n-1,n)$ with $n>5$ then $\cal N$ has type $\{n,10\}_{n-2}$ and is therefore in class~I, with ${\rm Aut}\,{\cal N}=A_n$. If $n=9$, for instance, the resulting class~III map $\cal M$ has type $\{630,630\}_{630}$ and characteristic $-1483791586099200$.

\section{Maps of small genus in class~III}

Each of the class~III maps constructed in the preceding section has rather large genus, whereas it interesting to know the lowest genera which can arise. Conder's lists of regular maps~\cite{Con} give the type $\{p,q\}_r$ of each map, for characteristic $\chi=-1,\ldots, -200$. A map in class~III must have $p=q=r$, so one can search through these lists for maps of such types and try to determine whether or not they are in class~III. In the case of non-orientable maps, this leads to the following result:

\begin{thm}
The only non-orientable class~III regular maps of characteristic $\chi\geq -200$ are Wilson's dual pair $\cal W$ and $D({\cal W})$ with $\chi=-70$, and their double covers with $\chi=-196$ constructed in Example~4A.
\end{thm}

\noindent{\sl Proof.} In Conder's list of non-orientable regular maps there are $21$ entries of type $\{p,p\}_p$, each representing a self-dual map or a dual pair. The information provided includes the multiplicities $mV$ and $mF$ of edges connecting adjacent pairs of vertices and adjacent pairs of faces. In a class~III map these must be equal, since the map is invariant under an operation sending vertices to faces. One can therefore exclude any of these $21$ entries with $mV\neq mF$, namely N22.3, N86.15, N170.15, N182.10 and N200.22. Conder's list also tells us whether or not an entry with $p=q$ represents a self-dual map: such a map cannot be in class~III, and this further criterion excludes the entries N12.3, N35.3, N44.6, N50.8, N119.6, N119.7, N146.8, N162.9, N200.23 and N200.26. Each entry gives defining relations for the automorphism group $G={\rm Aut}\,{\cal M}$, in terms of generators $R=r_0r_1$, $S=r_1r_2$ and $T=r_1$. Adding the relations $[R,S]=[R,T]=[S,T]=1$ gives a presentation for the abelianisation $G^{\rm ab}=G/G'$, and hence a set of normal generators for $\theta^{-1}(G')$ where $\theta:\Gamma\to G$ is the epimorphism $R_i\mapsto r_i$ corresponding to $\cal M$. Since $G'$ is a characteristic subgroup of $G$, if $\cal M$ is in class~III then $\theta^{-1}(G')$ must be a $\Sigma^+$-invariant subgroup of $\Gamma$, containing $\Gamma'$. As shown in \S 2, the only such subgroups are $\Gamma$, $\Gamma_{02}$, $\Gamma^*$ and $\Gamma'$, so if $\theta^{-1}(G')$ is not one of these then this entry can be excluded. This eliminates N62.4, N166.14 and N170.16. The remaining three entries are N72.9 of type $\{9,9\}_9$ corresponding to Wilson's class~III maps $\cal W$ and $D({\cal W})$, N198.6 of type $\{18,18\}_{18}$ corresponding to their class~III double covers described in Example~4A, and N119.5 of type $\{7,7\}_7$ and characteristic $-117$, with automorphism group $L_2(13)$. This last dual pair are in fact in class~II, having the self-dual map N119.6 as their third direct derivate; there is also another map N119.7 of type $\{7,7\}_7$ with automorphism group $L_2(13)$, but this is in class~IV.

To verify the assertions in this last sentence one needs to construct the epimorphisms $\Gamma\to G=L_2(13)$ corresponding to regular maps of type $\{7,7\}_7$. By applying a suitable automorphism of $G$ one may assume that
\[r_0=\pm\left(\begin{array}{cc}5&0\\ 0&-5\end{array}\right)
\quad{\rm and}\quad
r_2=\pm\left(\begin{array}{cc}0&1\\ -1&0\end{array}\right),
\quad{\rm while}\quad
r_1=\pm\left(\begin{array}{cc}a&b\\ c&-a\end{array}\right)\]
for some $a, b, c\in F_{13}$ with $a^2+bc=-1$. There are three conjugacy classes of elements of order $7$ in $G$, all self-inverse, consisting of the elements with traces $\pm 3$, $\pm 5$ and $\pm 6$. By applying this condition to the traces $\pm 3a$, $\pm(b-c)$ and $\pm 5(b+c)$ of the elements $r_0r_1$, $r_1r_2$ and $r_0r_1r_2$, and solving the resulting equations, one finds just the following solutions (unique up to multiplication by $-1$ and conjugation by an element of the Klein four-group $\langle r_0, r_2\rangle$):

\begin{enumerate}
\item $a=6$, $b=4$, $c=-6$, with $r_0r_1$, $r_1r_2$ and $r_0r_1r_2$ having traces $\pm 5$, $\pm 3$ and $\pm 3$;
\item $a=1$, $b=5$, $c=-3$, with $r_0r_1$, $r_1r_2$ and $r_0r_1r_2$ having traces $\pm 3$, $\pm 5$ and $\pm 3$;
\item $a=1$, $b=1$, $c=-2$, with $r_0r_1$, $r_1r_2$ and $r_0r_1r_2$ having traces $\pm 3$, $\pm 3$ and $\pm 5$;
\item $a=2$, $b=1$, $c=-5$, with $r_0r_1$, $r_1r_2$ and $r_0r_1r_2$ all having trace $\pm 6$.
\end{enumerate}

\noindent The dihedral group $\langle r_0, r_1\rangle$ of order $14$ is maximal in $G$, and it does not contain a Klein four-group, so in each case the elements $r_i$ generate $G$. It follows that there are exactly four regular maps ${\cal M}_1,\ldots, {\cal M}_4$ of type $\{7,7\}_7$ with automorphism group $G$, corresponding to these four solutions. They are non-orientable since $G$ has no subgroups of index $2$. Since ${\rm Aut}\,G=PGL_2(13)$ preserves traces, ${\cal M}_1$ and ${\cal M}_2$ form a dual pair, which must correspond to Conder's entry N119.5, while ${\cal M}_3$ and ${\cal M}_4$ are self-dual; similarly, ${\cal M}_3=P({\cal M}_1)$, so ${\cal M}_1$, ${\cal M}_2$ and ${\cal M}_3$ form a class~II orbit of $\Sigma$,  while ${\cal M}_4$ is in class~IV. This completes the proof, but one can go on to show that ${\cal M}_3$ is the map N119.6, so that ${\cal M}_4$ is N119.7, by verifying that in case~3 the generators satisfy the relation $(r_0r_1(r_2r_1)^2r_0r_1)^2=(RS^{-2}R)^2=1$ given in~\cite{Con} for N119.6, whereas in case~4 they do not. \hfill$\square$

\medskip

By contrast, the list of orientable regular maps in~\cite{Con} has $112$ entries with $p=q=r$, rather than the $21$ found above. Lemma~2.1 and its proof show that an orientable class~III map must have type $\{p,p\}_p$ with $p$ even, and that it must satisfy $G^{\rm ab}\cong V_8$; these and the preceding criteria eliminate all but R17.26, R33.45, R49.34, R61.13, R65.81, R73.85, R82.21, R97.108, R97.121 and R97.126 as possibilities for class~III maps. It would be possible, but very tedious, to deal with these cases by hand, as in the proof of Theorem~6.1. However, a computer search by Conder~\cite{Con09} has eliminated them all, by showing that the least genus of any orientable class~III map is $193$, attained by a dual pair of type $\{16,16\}_{16}$ with an automorphism group of order $2048$. His search (which also confirms Theorem~6.1) shows that this is the smallest possible automorphism group of such a map, and that the second smallest is the group $L_2(8)\times V_8$ of order $4032$ corresponding to the dual pair of maps of type $\{18,18\}_{18}$ and genus $392$ constructed in Example~4A.

It is interesting that in Conder's example the automorphism group is solvable (in fact nilpotent, having order $2^{11}$), whereas the examples constructed in this paper all have non-solvable automorphism groups. The methods of \S 4 yield further class~III maps with solvable groups, as coverings of Conder's, but it would be interesting to find infinite families of such maps which do not arise as coverings of other class~III maps. In particular one would like to see further examples with nilpotent groups (which are necessarily $2$-groups).

\end{document}